\documentclass{amsart}[12pt]
\usepackage[utf8]{inputenc}
\usepackage{amsmath}
\usepackage{amsfonts}
\usepackage{amssymb}

\newtheorem{theorem}{Theorem}[section]
\newtheorem{lemma}{Lemma}[section]
\newtheorem{corollary}{Corollary}[section]

\newtheorem{problem}{Problem}[section]

\newtheorem{remark}{Remark}[section]

\begin{document}
\title[Composition Operators ]%
{Composition Operators on Dirichlet Spaces over the Half-plane}

\author{Guangfu Cao, Haichou Li*}
\thanks{*Corresponding author}

\address{School of Mathematics and Information Science,
Guangzhou University
Guangzhou 510405,  China.}\email{ Email: guangfucao@163.com}

\address{College of Mathematics and informatics, South China Agricultural University, Guangzhou, 510640, China. }\email{ hcl2016@scau.edu.cn}

\begin{abstract}
As continuation of the study of polynomial approximation and composition operators on Dirichlet spaces of unit disk, which has settled a problem posed by Cima in 1976, the present paper aims to consider the case of the unbounded domains, such as the half-plane. Specifically, we may obtain the rational approximations in the Dirichlet spaces and characterize the composition operators which has dense range on the Dirichlet spaces over the half-plane. Moreover, this paper also considers the relationship between the Dirichlet spaces and Hardy spaces on half-plane.
\end{abstract}
\thanks{ Cao is supported by NSF of China (Grant No. 12071155); Li is supported by NSF of China (Grant No. 11901205). }

\maketitle

 \section { Introduction}

\ \ \ \
As we all know, there have been a great of studies on analytic composition operators acting on spaces of functions on the unit disc $\mathbb{D}$ of the complex plane $\mathbb{C}$, such as \cite{Cowen}, \cite{Shap}, \cite{Pau}, \cite{Zorb} and so on. However, composition operators on the equivalent spaces over the upper half-plane $\mathbb{C}^+$ enjoyed considerably less attention, except  \cite{Choe},  \cite{Elli},  \cite{Mat1}, \cite{Mat2}  and \cite{Shap2}. The composition operator $C_\varphi: f \rightarrow f \circ \varphi,$ where self-mapping $\varphi: \mathbb{C}^+ \rightarrow \mathbb{C}^+$ is analytic.  Although many corresponding spaces of analytic functions in the unit disc and half-plane are isomorphic, composition operators act very differently in this two setting. For example, the Littlewood principle shows that all composition operators acting on Hardy spaces of $\mathbb{D}$ are bounded, but Matache (\cite{Mat1}-\cite{Mat2}) proved that the Hardy spaces over $\mathbb{C}^+$ support no compact composition operators and not all are bounded. In addition, the Bergman spaces on $\mathbb{D}$ and $\mathbb{C}^+$ are also significantly different regarding the composition operators acting on them.

Based on the essential difference between the analytic composition operators acting on the spaces of functions on $\mathbb{D}$ and $\mathbb{C}^+$, this paper aims to consider the composition operators on the Dirichlet spaces $\mathfrak{D}(\mathbb{C}^+)$ over the upper half-plane.

In \cite{Cim}, J.A. Cima studied the Fredholmness of composition operators on  the Dirichlet spaces $\mathfrak{D}(\mathbb{D})$, and raised the following problem.

\begin{problem} \label{problem1} Characterize those $\varphi$ in the unit ball of $H^\infty$ for which $C_\varphi$ has dense range.
\end{problem}

Cima indicated that if $C_\varphi$ has dense range, then for $\psi(z)=r\varphi(z),\ 0<r<1,$ $C_\psi$  has dense range. He said that no conjecture on this question.  This is an old open problem relative  to the density of polynomials in the Dirichlet spaces on the domain $\Omega=\varphi(\mathbb{D})$. The first author of this paper and his coworkers (K.H. Zhu and L. He) find that the density of the polynomials in $\mathfrak{D}(\Omega)$ is equivalent to the density of the range of the composition operators $C_\varphi$ acting on $\mathfrak{D}(\mathbb{D})$, and obtain the following characterization results in \cite{Cao1}, which partly solves the open problem \ref{problem1}.\\

{\bf  Theorem A }\ (\cite{Cao1}) \ Suppose $\varphi: \mathbb{D} \rightarrow \mathbb{D}$ is analytic and $\Omega=\varphi(\mathbb{D})$. Then $C_\varphi$ is a bounded composition operator on Dirichlet spaces $\mathfrak{D}(\mathbb{D})$ with dense range if and only if the polynomials are dense in the Bergman space $A^2(\Omega)$ and the self-mapping $\varphi$ is univalent.\\

It is natural to ask how to solve this problem about the case of unbounded domain, more specifically, how to characterize the composition operators $C_\varphi$ with dense range on the Dirichlet spaces over the upper half-plane. This is a nontrivial problem, because there indeed exist some essential difficulties, such as finding dense set in $\mathfrak{D}(\mathbb{C}^+)$ is not so easy. Unlike $\mathfrak{D}(\mathbb{D})$ where the polynomials are dense, dense set in $\mathfrak{D}(\mathbb{C}^+)$ are often harder to come by. We first need some very smooth class of analytic functions that are dense in $\mathfrak{D}(\mathbb{C}^+)$ and will play the same role of the polynomials in the disc case $\mathfrak{D}(\mathbb{D})$, and then obtain the characterization of composition operators with dense range on $\mathfrak{D}(\mathbb{C}^+)$. Moreover, we also consider the boundedness  of composition operators acting on  $\mathfrak{D}(\mathbb{C}^+)$ and study the relationship between Hardy spaces $H^2(\mathbb{C}^+)$ and Dirichlet spaces $\mathfrak{D}(\mathbb{C}^+)$.

The present paper is arranged as follows. The second section is about some preliminary knowledge on the Dirichlet spaces. The third section consider the density of rational functions in $\mathfrak{D}(\mathbb{C}^+)$. The fourth section characterizes the composition operators with dense range on $\mathfrak{D}(\mathbb{C}^+)$. The last section considers the relationship between Dirichlet spaces $\mathfrak{D}(\mathbb{C}^+)$ and Hardy spaces $H^2(\mathbb{C}^+)$ over upper half-plane.

 \section { Preliminary Knowledge }

 Let  $\mathbb{D}$ be the open unit disc in the complex plane $\mathbb{C}$ and $\mathbb{C}^+$  the upper half-plane of complex plane, i.e.,
 $$\mathbb{D}=\{z\in \mathbb{C}: |z|<1\}$$
 and
 $$\mathbb{C}^+= \{z\in \mathbb{C}: \mathbf{Im}\ z>0\}<\infty.$$

 The Dirichlet spaces $\mathfrak{D}(\mathbb{C}^+)$ over the upper half-plane $\mathbb{C}^+$, denoted by $\mathfrak{D}(\mathbb{C}^+)$, consists of analytic functions $f$ on $\mathbb{C}^+$ satisfying
 $$\int_{\mathbb{C}^+}|f'(z)|^2d_{A(z)}$$
is finite.

 We here use $d_{A(z)}$ to denote Lebesgue (area) measure on $\mathbb{C}$ and its restriction to various domains, including $\mathbb{D}$ and $\mathbb{C}^+$.
 Moreover, we note that the norm of the function $f$ in Dirichlet spaces $\mathfrak{D}(\mathbb{C}^+)$ is defined as:
  $$\|f\|^2_{\mathfrak{D}(\mathbb{C}^+)}=|f(i)|^2+\int_{\mathbb{C}^+}|f'(z)|^2d_{A(z)}.$$
  It is easy to see that  $\mathfrak{D}(\mathbb{C}^+)$ is a Hilbert space whose inner production is given by
  $$\langle f, g \rangle=f(i)\overline{g(i)}+\int_{\mathbb{C}^+}f'(z)\overline{g'(z)}d_{A(z)}.$$

In order to state our main results, we still need to introduce the Dirichlet spaces and Bergman spaces of any domain in complex plane.
For any domain $\Omega\subseteq \mathbb{C}^+$, the Dirichlet space of $\Omega$, denote $\mathfrak{D}(\Omega)$, consists of all analytic functions $f$ on $\Omega$ satisfying that
$$\int_{\Omega}|f'(z)|^2d_{A(z)}$$
is finite.

The Bergman spaces $A^2(\Omega)$ on $\Omega$ consists all analytic functions  $f$ on $\Omega$ satisfying
$$\int_{\Omega}|f(z)|^2d_{A(z)}<\infty.$$
Let the fractional linear mapping

\begin{equation}\label{EQ1-1}
z=\alpha(w)=i\frac{1-w}{1+w}
\end{equation}
be the conformal mapping of unit disc $\mathbb{D}$ on to the upper half plane $\mathbb{C}^+$.

This mapping has inverse
$$w=\alpha^{-1}(z)=\frac{i-z}{i+z}$$
which from $\mathbb{C}^+$ to $\mathbb{D}$. We note that $d\theta$ is the arc-length Lebesgue measure on the unit circle
 $$\partial \mathbb{D}=\{w\in\mathbb{C}:|w|=1\},$$ and $dx$ is the Lebesgue measure on the real line $\mathbb{R}$.

We know that the fractional linear mapping $\alpha(w)$ maps the unit circle $\partial \mathbb{D}$ on to the real line $\mathbb{R}$. In fact, for any $w=e^{i\theta}\in\partial \mathbb{D}$, we have
$$\alpha(e^{i\theta})= i\frac{1-e^{i\theta}}{1+e^{i\theta}}=\frac{\sin\theta}{1+\cos\theta}\in \mathbb{R}.$$
Therefore, if we denote
$$\alpha(e^{i\theta})=i\frac{1+e^{i\theta}}{1-e^{i\theta}}=x,$$
then we have
$$e^{i\theta}=\frac{x-i}{x+i}$$
and obtain the following relation between $d\theta$ and $dx$ like that
\begin{equation}\label{EQ6}
d\theta=\frac{2}{1+x^2}dx,\ \  or,\ \  dx=\frac{1}{1+\cos \theta}d\theta.
\end{equation}

Let $\Phi$ be the class of all analytic self-mappings of $\mathbb{C}^+$. Each self-mapping $\varphi$ of $\mathbb{C}^+$ (i.e. $\varphi\in\Phi$) induces a composition operator $C_\varphi$ defined by
$$C_\varphi f=f\circ \varphi$$
for functions $f$ analytic on $\mathbb{C}^+$. It is clear that $C_\varphi$ takes the spaces of analytic functions on $\mathbb{C}^+$ into itself.

 Moreover, for each self-mapping  $\varphi$ of upper half-plane, we obtain a relative self-mapping $\phi$ of unit disc $\mathbb{D}$ as follows:
$$\phi = \alpha^{-1}\circ\varphi\circ\alpha.$$
For self-mapping $\varphi: \mathbb{C}^+ \rightarrow \mathbb{C}^+$, if $\Omega = \varphi(\mathbb{C}^+)$, that is, $\Omega \subset \mathbb{C}^+$ and  $\varphi: \mathbb{C}^+ \rightarrow \ \Omega $, we define the norm of the Dirichlet spaces $\mathfrak{D}(\Omega)$ on $\Omega$ by

$$\|f\|^2_{\mathfrak{D}(\Omega)}=|f(\varphi(i))|^2+\int_{\Omega}|f'(z)|^2d_{A(z)}.$$
It is easy to see that $\mathfrak{D}(\Omega)$ is also a Hilbert space with inner product:
$$\langle f, g \rangle=f(\varphi(i))\overline{g(\varphi(i))}+\int_{\Omega}f'(z)\overline{g'(z)}d_{A(z)}.$$

It is not hard to prove that the mapping $Uf=f\circ\varphi$ from Dirichlet spaces $\mathfrak{D}(\Omega)$ to $\mathfrak{D}(\mathbb{C}^+)$ is a unitary operator.

In order to consider the relationship between the Dirichlet spaces and  Hardy spaces, it is necessary to introduce the basic definitions of Hardy spaces on the upper half plane.

For $0<p<\infty$, let $H^p=H^p(\mathbb{C}^+)$ denote the Hardy space of analytic functions $f$ on the upper half plane $\mathbb{C}^+$ satisfying that the quantity

$$\|f\|_{H^p}=\sup\limits_{y>0}\left(\int_{-\infty}^{\infty}|f(x+iy)|^p dx\right )^{\frac{1}{p}}$$
is finite.

Moreover, we still need some knowledge about subharmonic function and harmonic majorant. A subharmonic function on an open set $O$ in the plane is a function $v:O\rightarrow [-\infty, \infty]$ such that

$1.$ $f$ is upper semicontinous: $f(z_0)\geq \lim\limits_{z\rightarrow z_0}f(z)$ for $z_0\in O$;

$2.$ For each $z_0\in O$ there is $f(z_0)>0$ such that the disc $$O(z_0,r(z_0))=\{z:|z-z_0|<r(z_0)\}$$ is contained in $O$ and  for every $r<r(z_0)$, there is
$$f(z_0)\leq \frac{1}{\pi r^2}\int\int_{|z-z_0|<r(z_0)}f(z)dxdy.$$

It is well known that every harmonic function is subharmonic, and for each analytic function $f(z)$ on set $O$ , the function $\log |f(z)|$ is subharmonic function, moreover, the function $|f(z)|^p=e^{p\log |f(z)|}$ is also a subharmonic function on $O$.

The subharmonic function $f(z)$ on $O$ has a harmonic majorant if there is a harmonic function $F(z)$ such that $f(z)\leq F(z)$ throughout set $O$.

 \section {Density of Rational Functions in Dirichlet Spaces }

 \ \ \ \ As we all know, the polynomials  play an indispensable role in Dirichlet spaces $\mathfrak{D}(\mathbb{D})$, because the density of polynomials in $\mathfrak{D}(\mathbb{D})$ is very important. Therefore, finding dense sets in Dirichlet spaces $\mathfrak{D}(\mathbb{C}^+)$ over upper half-plane is often very useful. Especially, considering the density of rational functions in  $\mathfrak{D}(\mathbb{C}^+)$ is necessary. However, unlike the classical Dirichlet spaces $\mathfrak{D}(\mathbb{D})$ in which the polynomials are always dense, the dense sets in $\mathfrak{D}(\mathbb{C}^+)$ are not so easy to found. In this section, our main work is to find and construct some dense sets in $\mathfrak{D}(\mathbb{C}^+)$.

 In order to obtain the dense sets of $\mathfrak{D}(\mathbb{C}^+)$, and since the conformal mapping $\alpha$ defined in (\ref{EQ1-1}) connects the unite disc and upper half-plane, we first have the relationship between $\mathfrak{D}(\mathbb{C}^+)$ and  $\mathfrak{D}(\mathbb{D})$ as the following lemma \ref{lemma1}. We have to note that the relationship between Dirichlet spaces of unite disc and upper half-plane is different from any other spaces of analytic functions, such as classical Hardy spaces and Bergman spaces, more details about those difference will be discussed in the section 5.

\begin{lemma}\label{lemma1}
An analytic function $f$ on $\mathbb{C}^+$ belongs to $\mathfrak{D}(\mathbb{C}^+)$ if and only if $f\circ\alpha\in \mathfrak{D}(\mathbb{D})$, where $\alpha$ is defined as (\ref{EQ1-1}) .
\end{lemma}

\begin{proof}
``Only if"  \ For any $w\in\mathbb{D}$, there exists a $z\in\mathbb{C}^+$ such that $z=\alpha(w).$
We have
\begin{equation*}
 \begin{split}
 \int_{\mathbb{D}}|[f\circ \alpha (w)]'|^2d_{A(w)} & = \int_{\mathbb{D}}|f'( \alpha (w))|^2|\alpha'(w)|^2d_{A(w)}\\
  & =\int_{\mathbb{C}^+}|f'(z)|^2d_{A(z)}.
\end{split}
\end{equation*}
The condition $f\in\mathfrak{D}(\mathbb{C}^+)$ implies the above formula is finite.

 Thus,  $f\circ \alpha (w)\in \mathfrak{D}(\mathbb{D})$ holds.

Similarly, it is not hard to  prove the ``If part" as follows.

``If part:" \ For any $z\in\mathbb{C}^+$, there exists a  $w\in\mathbb{D}$ such that $z=\alpha(w).$
So, we have
\begin{equation*}
 \begin{split}
\int_{\mathbb{C}^+}|f'(z)|^2d_{A(z)} & = \int_{\mathbb{D}}|f'( \alpha (w))|^2|\alpha'(w)|^2d_{A(w)}\\
  & =\int_{\mathbb{D}}|[f\circ \alpha (w)]'|^2d_{A(w)}.
\end{split}
\end{equation*}
Therefore, together with the condition that  $f\circ\alpha\in \mathfrak{D}(\mathbb{D})$, we obtain  that $f\in\mathfrak{D}(\mathbb{C}^+)$. The proof of Lemma \ref{lemma1} is complete.
\end{proof}

Based on Lemma \ref{lemma1}, if let the set $\mathbf{R}$ denote all the rational functions in $\mathfrak{D}(\mathbb{C}^+)$, then we can find that the rational functions set $\mathbf{R}$ is dense in $\mathfrak{D}(\mathbb{C}^+)$ and plays the same role as the polynomials in $\mathfrak{D}(\mathbb{D})$, that is the following theorem. Notice that it is not sure all the rational functions in upper half-plane $\mathbb{C}^+$  is dense in $\mathfrak{D}(\mathbb{C}^+)$, but only for the rational functions which belong to $\mathfrak{D}(\mathbb{C}^+)$. For example, $z\bar{\in}\mathfrak{D}(\mathbb{C}^+)$ or $\frac{z^2}{z+1}\bar{\in}\mathfrak{D}(\mathbb{C}^+)$.
\begin{theorem}\label{theorem1}
The set $\mathbf{R}$ which contains all the rational functions in Dirichlet spaces $\mathfrak{D}(\mathbb{C}^+)$ is dense in $\mathfrak{D}(\mathbb{C}^+)$.
\end{theorem}

\begin{proof}
Suppose the set $\mathbf{R}$ is not dense in $\mathfrak{D}(\mathbb{C}^+)$. Then there exists a function $f_0\in\mathfrak{D}(\mathbb{C}^+)$, for any rational functions sequence $\{r_k\}$ in set $\mathbf{R}$, such that
\begin{equation}\label{EQ1}
\|r_k-f_0\|_{\mathfrak{D}(\mathbb{C}^+)}\nrightarrow 0,  \ \ \ \ \ (k\rightarrow \infty).
\end{equation}
By $f_0\in\mathfrak{D}(\mathbb{C}^+)$, and  Lemma \ref{lemma1}, we have  $f_0\circ\alpha\in \mathfrak{D}(\mathbb{D})$.

Since the polynomials are dense in $\mathfrak{D}(\mathbb{D})$, for function $f_0\circ\alpha(w)\in \mathfrak{D}(\mathbb{D})$, there exists a sequence $\{p_k\}$ of polynomials in $\mathbb{D}$, such that
\begin{equation}\label{EQ2}
\|p_k-f_o\circ \alpha\|_{\mathfrak{D}(\mathbb{D})}\rightarrow 0,  \ \ \ \ \ (k\rightarrow \infty).
\end{equation}
Moreover, for any $w\in\mathbb{D}$, there exists a $z\in\mathbb{C}^+$ such that $z=\alpha(w) $ and  $w=\alpha^{-1}(z).$
From \eqref{EQ2}, the following holds
\begin{equation}\label{EQ3}
\|p_k\circ \alpha^{-1}-f_0\|_{\mathfrak{D}(\mathbb{C}^+)}\rightarrow 0,  \ \ \ \ \ (k\rightarrow \infty).
\end{equation}
Note that the sequence $\{p_k\circ \alpha^{-1}(z)\} $  is also contained in  $\mathbf{R}$.

It is clearly, \eqref{EQ3} is a contradiction with \eqref{EQ1}. Therefore, the set $\mathbf{R}$  is dense in $\mathfrak{D}(\mathbb{C}^+).$
\end{proof}

\begin{remark}\label{remark1}
Although Theorem \ref{theorem1} shows that  the rational function set  $\mathbf{R}$  is dense in $\mathfrak{D}(\mathbb{C}^+)$ and also implies that the set $\mathbf{R}$ will play the same role in $\mathfrak{D}(\mathbb{C}^+)$ as the polynomials play in $\mathfrak{D}(\mathbb{D})$, there is still some difference between  $\mathfrak{D}(\mathbb{C}^+)$ and $\mathfrak{D}(\mathbb{D})$. Because the fact that the polynomials is dense in $\mathfrak{D}(\mathbb{D})$ means that all the polynomials are belong to $\mathfrak{D}(\mathbb{D})$ and dense in it, however, not all the rational functions are the functions in $\mathfrak{D}(\mathbb{C}^+)$! So what kind  of rational functions belong to $\mathfrak{D}(\mathbb{C}^+)$?
\end {remark}

This question will be interest. That is,  it is necessary  to  consider such an important question: What is the special form of the rational functions in the set  $\mathbf{R}$? We find that more than one set satisfies this conditions. Firstly, from the density of the polynomials in $\mathfrak{D}(\mathbb{D})$ and  the relationship between $\mathfrak{D}(\mathbb{C}^+)$ and $\mathfrak{D}(\mathbb{D})$ (Lemma \ref{lemma1}), as well as using Theorem \ref{theorem1}, we can obtain one dense set of $\mathfrak{D}(\mathbb{C}^+)$ as  follows.

\begin{corollary}\label{corollary2}
Suppose $\alpha$ is the mapping defined in  (\ref{EQ1-1}) from $\mathbb{D}$ to $\mathbb{C}^+$, if $\mathbf{R}_{\alpha}(C^{+})$  is  a rational functions set defined as
  $$\mathbf{R}_{\alpha}(C^{+})=\{P\circ \alpha^{-1}(z):  z\in \mathbb{C}^+, P \ {\rm is \ any \  polynomial \  in \ the \ unite \ disc} \mathbb{D} \}$$
 then $\mathbf{R}_\alpha(C^{+})$  is a dense  subset of  $\mathfrak{D}(\mathbb{C}^+)$.
\end{corollary}

\begin{proof}
For any rational function $r(z)= p(\alpha^{-1}(z))\in \mathbf{R}_{\alpha}(C^{+})$,  together the relationship of the mapping $\alpha(w)$ and its reverse $\alpha^{-1}(z)$, where $w\in\mathbb{D}$ and $z\in \mathbb{C}^+$, we have  $p(w)=r\circ \alpha(w)$ is a polynomial in he unite disc $\mathbb{D}$. By the fact that any polynomial in $\mathbb{D}$ is a function of  $\mathfrak{D}(\mathbb{D})$, together with Lemma \ref{lemma1}, we obtain that the rational function $r(z)= p\circ \alpha^{-1}(z)$ is a function of $\mathfrak{D}(\mathbb{C}^+)$, and then, the set $\mathbf{R}_{\alpha}(C^{+})$ is a subset of $\mathfrak{D}(\mathbb{C}^+)$. Therefore, Theorem \ref{theorem1} shows that $\mathbf{R}_{\alpha}(C^{+})$ is dense in $\mathfrak{D}(\mathbb{C}^+)$. Thus, the proof of Corollary \ref{corollary2} is complete.
\end{proof}

From Corollary \ref{corollary2}, we obtain that some rational functions on the upper half-plane such as $P(\alpha^{-1}(z))$ must be the functions  of Dirichlet spaces over the upper half-plane. For more general case,  what kind of  the rational function is belong to the Dirichlet space $\mathfrak{D}(\mathbb{C}^+)$?  The second author and her coworkers \cite{LDQ} have considered this type question for Hardy spaces over tube domains, and obtain that the analytic rational functions with $L^{p}$  non-tangent boundary limit functions are belong to relative Hardy spaces. How about the Dirichlet spaces? Similarly, it will be found that the analytic rational functions with its derivative is $L^{2}$  non-tangent boundary limit functions are functions of Dirichlet  spaces of upper half-plane.

 \section {Characterization of Composition Operators with Dense Range }

\ \ \ \ After obtaining some dense subset of Dirichlet spaces on upper half-plane in  the last section, in what following, we are now to sole the main problem \ref{problem1} for the case of the upper half-plane as the following \ref{theorem2}, which is one characterization of  composition operators of $\mathfrak{D}(\mathbb{C}^+)$  with dense range, and also a generalization of Theorem A which is about the classical case of unit disc.

\begin{theorem} \label{theorem2}
Let $\varphi$ be an analytic self-mapping of $\mathbb{C}^+$ and $\Omega=\varphi(\mathbb{C}^+)$. Then, $C_\varphi$ reduced by $\varphi$ is a bounded composition operator on the Dirichlet spaces  $\mathfrak{D}(\mathbb{C}^+)$ with dense range if and only if $\varphi$ is univalent and the rational function  set $\mathbf{R}_{\alpha}(C^{+})$ is dense in the  Dirichlet spaces  $\mathfrak{D}(\Omega)$, where $\mathbf{R}_{\alpha}(C^{+})$ is defined as in Corollary \ref{corollary2}.
\end{theorem}

\begin{proof}
``Only if" \ First, if $C_\varphi$ is a bounded composition operator on  $\mathfrak{D}(\mathbb{C}^+)$ with dense range, then it is easy to show that $\varphi$ is univalent. In fact, if  $\varphi$ is not univalent, that is, there are two different points $z_1,z_2\in \mathbb{C}^+$ such that $\varphi(z_1)=\varphi(z_2)$. Moreover, for $z_1,z_2\in \mathbb{C}^+$, there exist two different points $w_1\in\mathbb{D}$ and $w_2\in\mathbb{D}$ such that $w_1=\alpha^{-1}(z_1)$ and $ w_2=\alpha^{-1}(z_2)$. For any function $f(z)\in\mathfrak{D}(\mathbb{C}^+)$, we have $C_\varphi f(z_1)=f(\varphi(z_1))$ and $C_\varphi f(z_2)=f(\varphi(z_2))$. So that the following formula holds
\begin{equation}\label{EQ5-1}
C_\varphi f(z_1)=C_\varphi f(z_2).
\end{equation}
By $f(z)\in \mathfrak{D}(\mathbb{C}^+)$ and Lemma \ref{lemma1}, we have $f\circ \alpha(w)\in\mathfrak{D}(\mathbb{D})$.
 Together with formula \eqref{EQ5-1}, we have
\begin{equation}\label{EQ5}
C_\varphi (f\circ \alpha)(w_1)=C_\varphi (f\circ \alpha)(w_2).
\end{equation}
However,  since $C_\varphi$ has dense range in  $\mathfrak{D}(\mathbb{C}^+)$, for any  function $f(z)\in\mathfrak{D}(\mathbb{C}^+)$, there exist a sequence $\{f_n(z)\}\subset \mathfrak{D}(\mathbb{C}^+)$ such that
$$\|C_\varphi f_n-f\|_{\mathfrak{D}(\mathbb{C}^+)}\rightarrow 0,\ \ \ \ \  (n\rightarrow\infty). $$
By Lemma \ref{lemma1}, there are $\{f_n\circ \alpha (w)\}\subset \mathfrak{D}(\mathbb{D})$ and  $f\circ \alpha (w)\in\mathfrak{D}(\mathbb{D})$. Thus, there is
$$\|C_\varphi (f_n\circ \alpha)-f\circ \alpha\|_{\mathfrak{D}(\mathbb{D})}\rightarrow 0,\ \ \ \ \  (n\rightarrow\infty).$$
If we take $f\circ \alpha (w)=w$, then $w_1=w_2$, moreover, there hold
$$z_1=\alpha(w_1)=\alpha(w_2)=z_2,$$
 which is contradiction with the condition $z_1\neq z_2$.
So that the mapping $\varphi$ is univalent.

 Next we are to prove the set $\mathbf{R}_{\alpha}(C^{+})$ is dense in $\mathfrak{D}(\Omega)$. If not so, then for any sequence $\{r_k\}$ of rational functions in $\mathbf{R}_{\alpha}(C^{+})$, there exists a function $g_0\in\mathfrak{D}(\Omega)$ such that,
\begin{equation}\label{EQ6-2}
\|r_k-g_0\|_{\mathfrak{D}(\Omega)}\nrightarrow 0,\ \ \ \ \  (k\rightarrow\infty).
\end{equation}
Thus, there is

\begin{equation}\label{EQ6-3}
\|C_\varphi r_k-C_\varphi g_0\|_{\mathfrak{D}(\mathbb{C}^+)}\nrightarrow 0,\ \ \ \ \  (k\rightarrow\infty).
\end{equation}
If let $f_0=C_\varphi g_0$, since $C_\varphi$ has dense range in $\mathfrak{D}(\mathbb{C}^+)$, then  there is a sequence $\{f_k(z)\}$ of functions in $\mathfrak{D}(\mathbb{C}^+)$ such that
\begin{equation}\label{EQ6-2-2}
\|C_\varphi f_k-f_0\|_{\mathfrak{D}(\mathbb{C}^+)}\rightarrow 0,\ \ \ \ \  (k\rightarrow\infty),
\end{equation}
where $f_k(z)=C_\varphi r_k=r_k\circ \varphi(z).$

By Theorem \ref{theorem1}, the set $\mathbf{R}_{\alpha}(C^{+})$ is dense in $\mathfrak{D}(\mathbb{C}^+)$, so that, for any function $f\in\mathfrak{D}(\mathbb{C}^+)$, there exist a sequence in $\mathbf{R}_{\alpha}(C^{+})$ converges to $f$ in   $\mathfrak{D}(\mathbb{C}^+)$ norm. Therefore, by $f_0\in \mathfrak{D}(\mathbb{C}^+)$, there is a sequence $\{r_k\}$ in $\mathbf{R}_{\alpha}(C^{+})$ such that
$$\|C_\varphi r_k-f_0\|_{\mathfrak{D}(\mathbb{C}^+)}\rightarrow 0,\ \ \ \ \  (k\rightarrow\infty).$$

However, the set $\mathbf{R}_{\alpha}(C^{+})$ is contained in Dirichlet spaces $\mathfrak{D}(\mathbb{C}^+)$, and then the rational function sequence $\{r_k\}$ is also in $\mathfrak{D}(\mathbb{C}^+)$. Obviously, this is  a contradiction with (\ref{EQ6-3}) and (\ref{EQ6-2-2}). Therefore, the rational function set $\mathbf{R}_{\alpha}(C^{+})$ is dense in $\mathfrak{D}(\Omega)$.

``If part " If $\varphi$ is univalent and $\mathbf{R}_{\alpha}(C^{+})$ is dense in $\mathfrak{D}(\Omega)$, we are to prove that the range, $R(C_\varphi)$, of composition operator $C_\varphi$, is dense in $\mathfrak{D}(\mathbb{C}^+)$.

 First the condition that $\varphi$ is univalent shows that $C_\varphi$ is univalent operator from $\mathfrak{D}(\Omega)$ to $\mathfrak{D}(\mathbb{C}^+)$ with the inverse $C_\varphi^{-1}$ from  $\mathfrak{D}(\mathbb{C}^+)$ to $\mathfrak{D}(\Omega)$. That is, for any function $g\in \mathfrak{D}(\mathbb{C}^+)$, there exists a function $f\in\mathfrak{D}(\Omega)$ such that $C_\varphi f=g.$

 Since $\mathbf{R}_{\alpha}(C^{+})$ is dense in $\mathfrak{D}(\Omega)$, there exist a rational function sequence $\{r_k\}$ contained in $\mathbf{R}_{\alpha}(C^{+})$ such that
 $$\|r_k- f\|_{\mathfrak{D}(\Omega)}\rightarrow 0, \ \ \ \ \ (k\rightarrow \infty).$$
 And then
  $$\|C_\varphi r_k- C_\varphi f\|_{\mathfrak{D}(\mathbb{C}^+)}\rightarrow 0, \ \ \ \ \ (k\rightarrow \infty).$$
Together with the relationship $C_\varphi f=g$, we have
 $$\|C_\varphi r_k- g|_{\mathfrak{D}(\mathbb{C}^+)}\rightarrow 0, \ \ \ \ \ (k\rightarrow \infty).$$
 Moreover, note that  $\mathbf{R}_{\alpha}(C^{+})$ is also contained in $\mathfrak{D}(\mathbb{C}^+)$, thus $\{r_k\}$ is also contained in $\mathfrak{D}(\mathbb{C}^+)$, hence, $R(C_\varphi)$ is dense in $\mathfrak{D}(\mathbb{C}^+)$. Therefore, the proof of Theorem \ref{theorem2} is complete.
\end{proof}

It should be pointed out that there is some difference between  Theorem A and Theorem \ref{theorem2}, separately about cases of the unit disc and the upper half-plane. For example, Theorem A can associate with Dirichlet spaces $\mathfrak{D}(\mathbb{D})$ and Bergmant spaces $A^2(\Omega)$, but Theorem \ref{theorem2} can not achieve that. This shows that the properties of composition operators on $\mathfrak{D}(\mathbb{D})$ are richer than $\mathfrak{D}(\mathbb{C}^+)$.

It is not hard to explain this difference. In fact, for the case of unit disc, as we all know, the polynomials play a very important role and they are both belong to the spaces $\mathfrak{D}(\mathbb{D})$ and $A^2(\mathbb{D})$. Moreover, since the derivative of any polynomial is till a polynomial, the polynomials are dense in $A^2(\Omega)$ if and only if they are dense in $\mathfrak{D}(\Omega)$, where $\Omega\subseteq \mathbb{D}.$ Therefore, Thoerem A characterizes the composition operators of $\mathfrak{D}(\mathbb{D})$  by using the  relation ship between  $\mathfrak{D}(\mathbb{D})$ and $A^2(\Omega)$. However, the case of half-plane couldn't. Because, for any rational functions $r(z)$ of $\mathfrak{D}(\mathbb{C}^+)$ on upper half-plane, we are not sure that $r(z)$ is also a  function of  Bergman spaces $A^2(\mathbb{C}^+)$. Therefore, under the same conditions of Theorem \ref{theorem2}, the density of the set $\mathbf{R}_{\alpha}(C^{+})$ in Bergman spaces $A^2(\Omega)$ could not hold.

Next, we continue to consider the boundedness of composition operators on $\mathfrak{D}(\mathbb{C}^+)$, and obtain the following result. Let's first recall a simple but important property of composition operators on Dirichlet spaces as the following Lemma \ref{lemma4.1}.

\begin{lemma}\label{lemma4.1} (\cite{Primer})
Let $D_1$ and $D_2$ be domains in complex plane. If $\varphi: D_1 \rightarrow D_2$ is univalent  and $f:D_2\rightarrow \mathbb{C}$ is a holomorphic function, then
$$\int_{D_1}|(f\circ \varphi)'(z)|^2d_{A(z)} =\int_{D_2}|f'(w)|^2d_{A(w)}.$$

\end{lemma}
In particular, taking $D_1=\mathbb{C}^+$ and $D_2=\varphi(\mathbb{C}^+)=\Omega$, we can obtain the other form of the norm of $\mathfrak{D}(\Omega)$ as
$$\|C_\varphi f\|^2_{\mathfrak{D}(\Omega)}= |f(\varphi(i))|^2+\int_{\Omega}|f'(w)|^2d_{A(w)}=|f(\varphi(i))|^2+\int_{\mathbb{C}^+}|(f\circ\varphi)'(z)|^2d_{A(z)}.$$

Based on the last lemma, we obtain the following result.
\begin{theorem} \label{theorem4}
If the self-mapping $\varphi: \mathbb{C}^+ \rightarrow \mathbb{C}^+$ is univalent, then the composition operator $C_\varphi$ is self-mapping from $\mathfrak{D}(\mathbb{C}^+)$ to $ \mathfrak{D}(\mathbb{C}^+)$ and  bounded.
\end{theorem}
\begin{proof}
Let $\Omega=\varphi (\mathbb{C}^+)\subseteq \mathbb{C}^+$, for any function $f$ of $\mathfrak{D}(\mathbb{C}^+)$, and by Lemma \ref{lemma4.1}, then we have
\begin{equation*}
 \begin{split}
\|C_\varphi f\|^2_{\mathfrak{D}(\mathbb{C}^+)}
&= |f(\varphi(i))|^2+\int_{\mathbb{C}^+}|(f\circ\varphi)'(z)|^2d_{A(z)} \\
&= |f(\varphi(i))|^2+\int_{\Omega}|f'(w)|^2d_{A(w)} \\
&\leq|f(\varphi(i))|^2+\int_{\mathbb{C}^+}|f'(w)|^2d_{A(w)}.
\end{split}
\end{equation*}
Since $f\in \mathfrak{D}(\mathbb{C}^+)$, there is  $$\int_{\mathbb{C}^+}|f'(w)|^2d_{A(w)}<\infty.$$
Therefore, $$\|C_\varphi f\|^2_{\mathfrak{D}(\mathbb{C}^+)}<\infty,$$ that is,  $C_\varphi f$ maps $\mathfrak{D}(\mathbb{C}^+)$  back into itself.

Moreover, by  the closed graph theorem, we obtain that $C_\varphi f$ is bounded on $\mathfrak{D}(\mathbb{C}^+)$. Thus, the proof of Theorem \ref{theorem4} is complete.
\end{proof}

In what following, we are to obtain some other properties of composition operators on $\mathfrak{D}(\mathbb{C}^+)$.

\begin{theorem} \label{theorem5}
Let $\Omega\subseteq \mathbb{C}^+$ be a simply connect domain and the mapping $\varphi: \mathbb{C}^+ \rightarrow \Omega$ be univalent.  If the Lebesgue measure of the complement of $\Omega$ in $\mathbb{C}^+ $is zero ($m(\mathbb{C}^+\setminus \Omega)=0$), then the composition operator $C_\varphi: \mathfrak{D}(\mathbb{C}^+)\rightarrow \mathfrak{D}(\mathbb{C}^+)$ is bounded and has closed range.
\end{theorem}

\begin{proof}
Since the mapping $\varphi$ is univalent, Theorem \ref{theorem4} shows that the composition operator $C_\varphi: \mathfrak{D}(\mathbb{C}^+)\rightarrow \mathfrak{D}(\mathbb{C}^+)$ is bounded.

Consider the norms of $\mathfrak{D}(\mathbb{C}^+)$ and $\mathfrak{D}(\Omega)$ respectively as follows:

$$\|f\|^2_{\mathfrak{D}(\mathbb{C}^+)}=|f(i)|^2+\int_{\mathbb{C}^+}|f'(z)|^2d_{A(z)}$$
and
$$\|f\|^2_{\mathfrak{D}(\Omega)}=|f(\varphi(i))|^2+\int_{\Omega}|f'(z)|^2d_{A(z)}.$$
By the condition $m(\mathbb{C}^+\setminus \Omega)=0$, we obtain the following formula for any function $f$ of $\mathfrak{D}(\mathbb{C}^+)$:
$$\int_{\mathfrak{D}(\mathbb{C}^+)}|f'(z)|^2d_{A(z)}=\int_{\Omega}|f'(z)|^2d_{A(z)}+\int_{\mathbb{C}^+ \setminus \Omega}|f'(z)|^2d_{A(z)}=\int_{\Omega}|f'(z)|^2d_{A(z)}.$$

Therefore, it is easy to prove that the two norms $\|f\|^2_{\mathfrak{D}(\mathbb{C}^+)}$ and $\|f\|^2_{\mathfrak{D}(\Omega)}$ are equivalent on $\mathfrak{D}(\mathbb{C}^+).$
Since the rational function set $\mathbf{R}_{\alpha}(C^{+})$ defined in Corollary \ref{corollary2} is contained in $\mathfrak{D}(\mathbb{C}^+)$, and $\mathfrak{D}(\mathbb{C}^+)$ is a subspace of $\mathfrak{D}(\Omega)$, we know that $\mathbf{R}_{\alpha}(C^{+})$ is also contained in $\mathfrak{D}(\Omega)$.

Moreover, since $\mathbf{R}_{\alpha}(C^{+})$ is dense in $\mathfrak{D}(\mathbb{C}^+)$, the closure of $\mathbf{R}_{\alpha}(C^{+})$ is just $\mathfrak{D}(\mathbb{C}^+),$  thus  $\mathfrak{D}(\mathbb{C}^+)$ is a closed subspace of $\mathfrak{D}(\Omega)$. We note that  the mapping $C_\varphi f=f\circ \varphi : \mathfrak{D}(\Omega)\rightarrow \mathfrak{D}(\mathbb{C}^+)$ is unitary operator. Therefore, the range of $C_\varphi$ is closed.
\end{proof}

 \section {The Relationship Between Dirichlet Spaces and Hardy Spaces on Half-plane }

\ \ \ \  From Lemma \ref{lemma1} of Section 3, we obtain the relationship between the Dirichlet spaces $\mathfrak{D}(\mathbb{D})$ and  $\mathfrak{D}(\mathbb{C}^+)$ by making essential use of the linear fractional map $\alpha(w)=i\frac{1-w}{1+w}$  which  transforms  unit disc on to upper-half plane.

There are some remarks we might make. The ``If and only if" relationship  between the Dirichlet spaces for the disc and half-plane in Lemma \ref{lemma1}, only holds for Dirichlet spaces but not for any other analytic function spaces, such as Hardy spaces and Bergman spaces. Of course, all the differences are due to their different definitions.

 For example, especially for classical Hardy spaces $H^p$, $0<p<\infty$, the relationship between Hardy spaces $H^p(\mathbb{D})$ and  $H^p(\mathbb{C}^+)$ is not so simple and beautiful as Dirichlet spaces in Lemma \ref{lemma1}. There are many discusses about the relation between the Hardy spaces for the disc and half-plane, such as Cima \cite{Cim2}, Garnett \cite{Gar} and Hoffman \cite{Hoffman}. In \cite{Cim2}, it shows that any analytic function $f$ of upper half-plane $\mathbb{C}^+$ belongs to Hardy spaces $H^p(\mathbb{C}^+)$ if and only if $(f\circ \alpha) (\alpha')^{1/p}\in H^p(\mathbb{D}),$ where $\alpha$ is defined as (\ref{EQ1-1}). In \cite{Gar}, only the ``only if part" as in Lemma \ref{lemma1} holds for Hardy spaces, that is, if $f(z)\in H^p(\mathbb{C}^+)$ then $f\circ \alpha(w) \in H^p(\mathbb{D})$, the converse is not true. We still take some information from reference \cite{Gar} and \cite{Hoffman} which imply that $f(w)\in H^p(\mathbb{D}$ if and only if $\frac{\pi^{-1/p}}{(z+i)^{2/p}}f(\alpha^{-1}(z))\in H^p(\mathbb{C}^+)$.

In this section, we mainly to consider the relationship between Dirichlet spaces $\mathfrak{D}(\mathbb{C}^+)$ and classical Hardy spaces $H^2(\mathbb{C}^+)$ of half-plane. It is well known, for the case of unit disc, the following relationship between $\mathfrak{D}(\mathbb{D})$ and  $H^{2}(\mathbb{D})$ holds.

\begin{lemma}\label{lemma2}
The  Dirichlet spaces $\mathfrak{D}(\mathbb{D})$ is contained in Hardy spaces $H^{2}(\mathbb{D})$.
\end{lemma}

It is nature to ask how about the relationship between Dirichlet spaces $\mathfrak{D}(\mathbb{C}^+)$ and Hardy spaces $H^{2}(\mathbb{C}^+)$ on the upper half-plane, is it the same as the case of unit disc? We find that the answer is ``no" and obtain the similar result such as Theorem \ref{theorem6}, which is simple but very useful for the theory of Dirichlet spaces $\mathfrak{D}(\mathbb{C}^+).$ Notice that the proof of Theorem \ref{theorem6} is  different from the proof of Lemma \ref{lemma2}, because the unbounded of upper half-plane restrict that we could not define the  Dirichlet spaces $\mathfrak{D}(\mathbb{C}^+)$ by the Taylor coefficients like $\mathfrak{D}(\mathbb{D}).$ Therefore, we have to find some other  methods to prove it, although it is not difficult. Moreover, in order to prove  Theorem \ref{theorem6}, we will make use of the previous stated relationship  between $H^{2}(\mathbb{D})$  and $H^{2}(\mathbb{C}^{+})$, which is the bridge between the case of unit disc and upper half-plane.

More precisely, we need the following lemmas.

\begin{lemma}\label{lemma3}(\cite{Gar})
If $f(z)$ is an analytic function in upper half-plane such that the subharmonic function $|f(z)|^p$, $0<p<\infty$ has a harmonic majorant, then
$$F(z)=\frac{\pi^{-1/p}}{(z+i)^{2/p}}f(z)$$
is in Hardy spaces $H^p(\mathbb{C}^+)$.
\end{lemma}

\begin{lemma}\label{lemma4}(\cite{Gar})
Let $g(w)$ be a subharmonic function in the unit disc. Then $g(w)$ has a harmonic majorant if and only if
$$\sup\limits_{0\leq r \leq 1}\int g(re^{i\theta})d\theta <\infty.$$
\end{lemma}

Clearly, Lemma \ref{lemma4} implies that for function $g(w)$ of Hardy spaces $H^p(\mathbb{D})$, $|g(w)|^p$ has a harmonic majorant.

 Based on the above  Lemma results, we are now to obtain the following relationship between the Hardy spaces $H^{2}(\mathbb{C}^+)$ and Dirichlet space $\mathfrak{D}(\mathbb{C}^+)$.
\begin{theorem} \label{theorem6}
If the analytic function  $f(z)\in\mathfrak{D}(\mathbb{C}^+)$, then  $f(z)[(\alpha^{-1}(z))']^{1/2}\in H^{2}(\mathbb{C}^+)$, that is, $\frac{\sqrt{-2i}}{z+i} f(z)\in H^{2}(\mathbb{C}^+)$, where $\alpha$ is defined as (\ref{EQ1-1}).
\end{theorem}

\begin{proof}
On one hand, for any function $f(z)\in \mathfrak{D}(\mathbb{C}^+)$,  it is easy to obtain that the functions $f(z)[(\alpha^{-1}(z))']^{1/2}$ and $\frac{\sqrt{-2i}}{z+i} f(\alpha^{-1}(z))$ are also analytic on upper half-plane $\mathbb{C}^+$.

On the other hand, by $f(z)\in \mathfrak{D}(\mathbb{C}^+)$ and  Lemma \ref{lemma1}, we have  $f\circ \alpha (w)\in \mathfrak{D}(\mathbb{D})$.  Consequently, Lemma \ref{lemma2} shows that  $f\circ \alpha (w)\in H^{2}(\mathbb{D})$.

Moreover, Lemma \ref{lemma4} induces  that $|f\circ \alpha(w)|^p$ has a harmonic majorant $u(w)$. Hence, $|f\circ \alpha(\alpha^{-1}(z))|^p=|f(z)|^p$ has a harmonic majoran $u(\alpha(z))$.

Therefore, by Lemma \ref{lemma3}, we obtain that $\frac{\pi^{-1/2}}{z+i}f(z)\in H^{2}(\mathbb{C}^+)$.  So that, there also holds   $\frac{\sqrt{-2i}}{z+i} f(z)\in H^{2}(\mathbb{C}^+).$

At last, we note that $(\alpha^{-1})'(z)=\frac{-2i}{(z+i)^2}$, then the relationship
 $$f(z)[(\alpha^{-1}(z))']^{1/2}\in H^{2}(\mathbb{C}^+)$$
is also established.
\end{proof}

%\par
%\vskip 0.5 cm
%{\bf Acknowledgement}\ \ The author wants to thank the
%referee for his careful reading of the manuscript and many useful
%comments for improvements of the presentation.

\vskip 0.5cm{

\end{document}